\documentclass{amsart}

\usepackage{amsmath, amssymb,epsfig, verbatim}
\usepackage{graphicx}
\usepackage{epic}
\usepackage{color}

\newtheorem{theorem}{Theorem}[section]

\newtheorem{lemma}[theorem]{Lemma}
\newtheorem{proposition}[theorem]{Proposition}

\newtheorem{remark}[theorem]{Remark}
\numberwithin{equation}{section}

\newcommand{\QQ}{\mathbf Q}
\newcommand{\ZZ}{\mathbf Z}
\newcommand{\RR}{\mathbf R}

\newcommand{\FF}{\mathbf F}
\newcommand{\HH}{\mathbf H}
\newcommand{\cB}{\mathcal B}

\newcommand{\nrd}{\operatorname{nrd}}
\newcommand{\trd}{\operatorname{trd}}
\newcommand{\vol}{\operatorname{vol}}

\newcommand{\PSL}{\operatorname{PSL}}
\newcommand{\SL}{\operatorname{SL}}
\newcommand{\Fix}{\operatorname{Fix}}

\newcommand{\sign}{\operatorname{sign}}

\newcommand{\ii}{\mathrm i}
\newcommand{\CP}{\mathbf{CP}}
\newcommand{\RP}{\mathbf{RP}}

\begin{document}

\title{Simply connected quotients of fake quadrics}
\author{Zolt\'an Szab\'o}
\address{Department of Mathematics, Princeton University,
Princeton, NJ 08544}
\email{szabo@math.princeton.edu}
\date{November 6, 2026}

\begin{abstract}
  Fake quadrics are minimal complex surfaces of general type with   $b_1=0$ and $e=4$.
  The paper gives a short introduction to their arithmetic construction, aimed at low-dimensional topologists,
  and studies their holomorphic and antiholomorphic quotients and their fundamental groups. The running example in the paper comes from a quaternion algebra over a totally real
cubic field. The lattice is generated by three explicit pairs of matrices.
A holomorphic involution gives (after resolving the singularities) a simply connected minimal complex surface of general type
  with $p_g=0$ and $K^2=4$, a smooth $4$-manifold homeomorphic but non-diffeomorphic to $\CP^2$ blown up $5$ times.
 An antiholomorphic involution gives a smooth four-manifold homeomorphic to $S^4$, with a knotted torus as branch locus.

\end{abstract}

\maketitle

\section{Introduction}

The complex surface $\CP^1\times \CP^1$ is called a quadric since it can be given by the equation
$x^2+y^2+z^2+w^2=0$ in $\CP^3$.
A fake quadric is a smooth minimal compact complex surface of general type
with $p_g=q=0$ and $K^2=8$. To explain the name, note that  these minimal surfaces have the same Euler characteristic, signature and rational homology
groups as $\CP^1\times \CP^1$. (Note $b_2^+=1+2p_g$ and $b_1=2q$.) 
A well-known construction of fake quadrics is to divide
$\HH\times\HH$ by a suitable torsion-free irreducible arithmetic
lattice, see \cite{Shavel,Takeuchi90,Dzambic,DR,LSV}. See also \cite{Cat2} for an overview on the reducible case where the surface admits a finite cover that is a product of two curves. See also \cite{Rito, BR}
for a more recent and different construction.

The goal of this paper is
to give a short introduction to the irreducible arithmetic construction  aimed at
low-dimensional topologists. The running example uses the cubic field $F= \QQ(t)$ with
$t^3-t^2-4t-1=0$ and a quaternion algebra over $F$. The fundamental group is the norm-one subgroup of a maximal order, divided by $\pm 1$.

We work with an explicit description
of one such surface by three pairs of matrices in
$\SL_2(\RR)\times\SL_2(\RR)$. The number field, quaternion algebra
and order behind this construction are explained in
Sections $2$ and $3$. Section $4$ uses a volume bound computation to prove
that the displayed matrices generate the required lattice.

Let $t_0<t_1<t_2$ be the three real roots of $t^3-t^2-4t-1=0$ and

$$ a(s)=3s^2-4s-12,\ \  b(s)=-2s^2+4s+4. $$

At $s=t_1,t_2$ we have $a(s)<0<b(s)$ and at $s=t_0$ both $a(s)$ and
$b(s)$ are negative. For $s=t_1,t_2$ let
$$ \alpha=\sqrt{-a(s)},\qquad \beta=\sqrt{b(s)}.$$

We will use the following three matrices:
\[
A_s=
\begin{pmatrix}
 -s-1 & \alpha+(s^2-2)\beta\\
 -\alpha+(s^2-2)\beta & -s-1
\end{pmatrix},
\]
\[
B_s=\frac12
\begin{pmatrix}
 -s^2 & (s+2)\alpha+(3s+4)\beta\\
 -(s+2)\alpha+(3s+4)\beta & -s^2
\end{pmatrix},
\]
\[
C_s=\frac14
\begin{pmatrix}
 2s^2-4s-8-\alpha\beta & -s\beta\\
 -s\beta & 2s^2-4s-8+\alpha\beta
\end{pmatrix}.
\]

It is easy to check that these matrices all have determinant one. We get three pairs
$$ A =(A_{t_1}, A_{t_2}) \ \ B=(B_{t_1},B_{t_2}) \ \ C= (C_{t_1}, C_{t_2}) $$
in $\SL_2(\RR)\times \SL_2(\RR)$. Let $\Gamma$ denote the group generated by their projective images in
$\PSL_2(\RR)\times \PSL_2(\RR)$. Elements in $\Gamma$ act on $\HH \times \HH$ by the usual formulas:
$$(z_1,z_2)\mapsto ((a_1z_1+b_1)/(c_1z_1+d_1),(a_2z_2+b_2)/(c_2z_2+d_2)).$$
In $\HH$ we use the hyperbolic metric with curvature $-1$. Then these actions on $\HH\times \HH$ are isometries.
The volume of the fundamental domain is related to the (orbifold) Euler characteristic of the quotient.

Let $X=\Gamma\backslash(\HH\times\HH)$. We will use the following standard
facts, see \cite{Shavel} and \cite[Section~2.1]{Dzambic}.

\begin{proposition}\label{lattice}
Let $\Gamma\subset\PSL_2(\RR)\times\PSL_2(\RR)$ be an irreducible, cocompact,
torsion-free lattice. Then $X=\Gamma\backslash(\HH\times\HH)$ is a minimal
surface of general type with $q(X)=0$ and $e(X)=\vol(X)/4\pi^2$. In
particular, $X$ is a fake quadric if and only if $\vol(X)=16\pi^2$.
\end{proposition}

In Section $4$ we will prove the following:

\begin{theorem}\label{thm1}
The group $\Gamma$ is a cocompact, torsion-free, irreducible lattice in 
$\PSL_2(\RR)\times \PSL_2(\RR)$. Furthermore the volume of $\Gamma\backslash(\HH\times \HH)$ is equal to $16\pi^2$.
\end{theorem}

In Sections $3$ and $4$ we will show that  $\Gamma$ is the group of norm-one elements, modulo $\pm 1$ of a maximal order in a quaternion algebra
over the cubic field $\QQ(t)$, $t^3-t^2-4t-1=0$. The commensurability class is listed in \cite{LSV} and given as $(n,d,D,N)=(3,169,13,1)$.
Here $n$ is the degree of the field, $d$ is the discriminant of the field, $D$ is the norm of the discriminant of the quaternion algebra, and
$N=1$ means that it is a maximal order.

Proposition~\ref{lattice} and a short computation
 using some of the relations among $A,B,C$ give:

\begin{theorem}\label{thm2}
  The quotient    $X =\Gamma\backslash(\HH\times \HH)$ is a minimal complex surface of general type with
  $$H_1(X;\ZZ)= \ZZ/7, \ \ e(X)=4.$$
\end{theorem}

  Dividing by finite group actions (and then desingularizing) is a standard tool to construct new complex
  surfaces from old ones. In the context of fake quadrics see
for example \cite{DR}, see also \cite{Rito}, \cite{BR}.
Using the same idea for $X$ one gets:

\begin{theorem}\label{thm3}
  There is a simply connected  minimal complex surface of general type $Z$ with
  $p_g(Z)=0$ and $K^2=4$ that contains $4$ disjoint holomorphic spheres with self-intersection $-2$. The surface $Z$ is the minimal resolution of the
  quotient of $X$ by a holomorphic involution. Furthermore $Z$ is homeomorphic but not diffeomorphic to $\CP^2\#5\overline{\CP}^2$.
  
\end{theorem}

Of course the last part is standard.
First note that $p_g=0$ is equivalent to $b_2^+=1$, and then $4=K^2=2e+3\sign$ and $\pi_1=1$ show  that $Z$ is homotopy equivalent to
$\CP^2\#5\overline{\CP}^2$. Freedman's theorem proves that they are also homeomorphic. 
Since a minimal surface of general type  has
nonzero Seiberg-Witten invariants, see for example \cite{FM}, $Z$ is not diffeomorphic to
$\CP^2\#5\overline{\CP}^2$.
 
Simply connected minimal surfaces of general type with $p_g=0$ and $K^2=4$ were already known,
see \cite{PPS4}, obtained by $\QQ$-Gorenstein smoothing \cite{LP}.
The surface $Z$ is not biholomorphic to the examples of \cite{PPS4}: those surfaces have 
ample canonical class  showing that they don't  contain $(-2)$-curves. On the other hand $Z$
contains four such disjoint holomorphic $-2$ curves.
See also \cite{BCP} for a survey of complex surfaces of
general type with $p_g=0$.

\begin{remark}
  In order to keep  this paper short we concentrate on  a single example. 
There are finitely many fake quadrics that come from the irreducible arithmetic lattices. Some of them are described in \cite{Shavel,Takeuchi90,Dzambic,DR,LSV},
and \cite{LSV} gives a list of the quaternion algebras that may contain additional examples.
Some of the arithmetic fake quadrics  admit holomorphic
involutions, see \cite{DR}  and the approach given in Sections $4$ and $5$ can  be used to study the fundamental groups or the first homology groups
of the quotients. In lots of cases these groups are non-trivial.

There are, however, further examples whose holomorphic involution quotients have
simply connected minimal resolutions of general type  with $p_g=0$ and $K^2=4$. Just as for $Z$, these surfaces are homeomorphic
but not diffeomorphic to $\CP^2\#5\overline{\CP}^2$, and they contain
four disjoint $(-2)$-curves. Whether they are diffeomorphic to one another or to $Z$, and how their deformation types compare,
are interesting questions.

Fake quadrics with larger holomorphic symmetry groups provide
another direction, see \cite{DR}. Suitable actions of $\ZZ/2\times\ZZ/2$ and $\ZZ/2\times\ZZ/4$ give simply connected minimal
quotient resolutions of general type
with $p_g=0$ and $K^2=2$ and $K^2=1$, respectively. These
give exotic smooth structures on $\CP^2$ blown up seven or eight times.
We plan to discuss these further examples and their fundamental group calculations in a follow-up paper.

 \end{remark}

Another construction that is also interesting in smooth $4$-manifold topology is to divide a complex surface
by a real involution. An antiholomorphic involution $\tau$ on a complex surface $S$
gives a pair $(Y,\Sigma)$: the smooth
four-manifold $Y=S/\langle\tau\rangle$ and the branch locus $\Sigma$. Sometimes these are standard. For example, complex
conjugation on $\CP^2$ has quotient $S^4$ with an
embedded $\RP^2$ as branch locus \cite{Arnold}. For $\CP^1\times\CP^1$, the quotient
by $\tau(z_1,z_2)=(\overline{z}_1,\overline{z}_2)$ is again $S^4$, with a standardly embedded torus as branch locus.

In other cases the embedded surface is exotic. For example, \cite{FKV2} gives exotic embeddings
of the connected sum of ten copies of $\RP^2$ in $S^4$ as
branch loci of real quotients of Dolgachev surfaces. See also \cite{Finashin}
and \cite{MORSU} for related constructions in the symplectic category, and \cite{SS} for the
antiholomorphic involution on the
Cartwright-Steger surface.

\begin{theorem}\label{thm4}
  The fake quadric $X$ admits an antiholomorphic involution $\tau$. The quotient $Y=X/\tau$ is a smooth $4$-manifold homeomorphic
  to the four-dimensional sphere $S^4$, and the branch locus is a knotted torus.
\end{theorem}

{\it Motivation.} This work started as a study of antiholomorphic involutions
on fake quadrics, along the lines of \cite{SS} where we study the Cartwright-Steger surface.
The aim was to understand the branch locus $\Sigma$ in the real quotients $Y$, and
to find interesting pairs $(Y,\Sigma)$ with $Y$  homeomorphic to $S^4$. 
As it turns out, the same method to compute the fundamental group works for  holomorphic involutions as well. This then led to  
Theorem~\ref{thm3}. 

{\it Acknowledgments.}
The author would like to thank Andr\'as Stipsicz for helpful discussions. This work was
partially supported by the Simons Collaboration Grant on New
Structures in Low-Dimensional Topology.
Part of this paper was prepared as background material for the {\it Exotica II Workshop}
at Caltech on  November 12-16, 2026, supported by the same collaboration.

{\it Statement on the use of AI.}
There are programs attached
to the arXiv version of this paper. These are used for
giving an upper bound for the volume of the fundamental domain for the lattice used in the paper, and for checking some relations
in this lattice. The AI assistants ChatGPT and Claude were used to
write these programs and do preliminary calculations for some other quaternion algebras listed in \cite{LSV}.
 These LLMs  were also used for typesetting and language
editing. The author assumes intellectual responsibility for the
mathematical statements and for any remaining errors.

\section{Arithmetic lattices}

This section is an overview of the arithmetic construction behind
Section $3$. See \cite{Shavel,Dzambic,DR,LSV} for this construction and
\cite{Voight} for background on quaternion algebras.

Let $f(T)$ be an irreducible monic integer polynomial with only real
roots. Then $F=\QQ[T]/(f)$ is a totally real number field. Write $t$ for the image
of $T$ and $R$ for its ring of algebraic integers,
the elements satisfying a monic polynomial over $\ZZ$. Each root gives an
embedding $F\hookrightarrow\RR$ by substitution. In the cubic example of Section~3, we use
$f(T)=T^3-T^2-4T-1$. In this case the discriminant is $169$, and
$f(T-4)=T^3-13T^2+52T-65$
satisfies the Eisenstein condition at $13$, so $R=\ZZ[t]$. 

For nonzero $a,b\in R$, the quaternion algebra $\cB=(a,b)_F$ has
basis $1,i,j,ij$ over $F$, with
$$ i^2=a,\qquad j^2=b,\ \  ij=-ji,\ \ (ij)^2=-ab $$
We also require that $\cB$ is a division algebra: every nonzero element
has an inverse.

Note that in the construction,  $a,b \in F-\{0 \}$.  However multiplying $a$, $b$ by squares of nonzero elements leads to
isomorphic quaternion algebras, so we can assume $a,b \in R$. After fixing $a,b$ we still need to check that $\cB$ is a division algebra.
In practice we will work over the reals or  modulo a prime.

It is easy to check at the real roots of $f$. Given such a root, we get a 
 real embedding $\sigma$. Now if both $\sigma(a),\sigma(b)$
are negative then we can rescale and see that the scalar extension 
is isomorphic to the standard  Hamiltonian quaternions.
 Otherwise it is isomorphic to the algebra of two by two matrices $M_2(\RR)$.
We say that $\cB$ is ramified at $\sigma$ in the first case
and split at $\sigma$ in the second. 

We can also work over prime ideals in $R$.
For example $f(T)=T^3-T^2-4T-1$ has a root modulo $13$ so we can take the map $R\rightarrow \FF_{13}$ given by $t\mapsto -4$ since $f(-4)\equiv0\pmod{13}$.
Then the kernel is a prime ideal
$\mathfrak{p}$. Let $F_{\mathfrak p}$ denote the completion of $F$ at $\mathfrak p$. The quaternion algebra
$\cB$ is called split at $\mathfrak p$ if $\cB\otimes _F F_{\mathfrak p}$ is isomorphic to the matrix algebra
$M_2(F_{\mathfrak p})$. Otherwise it is called ramified at $\mathfrak p$.

We require exactly two split real embeddings, these are denoted by
$\sigma_1,\sigma_2$. In degree $3$ or more that immediately implies that $\cB$ is a division algebra.
Different choices of $a,b$ could still lead to
isomorphic division algebras. What really matters is the list of primes and real places over which $\cB$ ramifies. Note that
this list is always finite and even, see \cite{Voight}.

In our example $\sigma_n(a)<0<\sigma_n(b)$ for $n=1,2$. At
either embedding put $\alpha=\sqrt{-\sigma_n(a)}$ and
$\beta=\sqrt{\sigma_n(b)}$, as in the introduction. We use the
matrices

\begin{equation}
  \label{matrix1}
i\mapsto \begin{pmatrix}0&-\alpha\\\alpha&0\end{pmatrix},\ \ j\mapsto \begin{pmatrix}0&-\beta\\-\beta&0\end{pmatrix},\ \ 
ij\mapsto \begin{pmatrix}\alpha\beta&0\\0&-\alpha\beta\end{pmatrix}.
\end{equation}

Extended $\sigma_n$-linearly, \eqref{matrix1} defines an injective $F$-algebra
homomorphism $\cB\to M_2(\RR)$ and an isomorphism
$\cB\otimes_{F,\sigma_n}\RR\cong M_2(\RR)$, see \cite[Section~2.3]{Voight};
in particular $O^1\to\PSL_2(\RR)\times\PSL_2(\RR)$ has kernel $\{\pm1\}$.

For $q=x_0+x_1i+x_2j+x_3ij$, its conjugate is
$\bar q=x_0-x_1i-x_2j-x_3ij$. The reduced trace and norm of $q$ are
$$ \trd(q)=2x_0,\ \ \nrd(q)=q\bar q=x_0^2-ax_1^2-bx_2^2+abx_3^2.$$

They become matrix trace and determinant at a split embedding, and $q^2-\trd(q)q+\nrd(q)=0$.

A real invertible matrix with positive determinant acts on the upper half-plane by the usual  $z\mapsto (uz+v)/(wz+r)$. In case the determinant is negative
$z\mapsto (u{\overline z}+v)/(w{\overline z}+r)$. A nonzero quaternion element $q$ acts by its two matrices. It acts holomorphically when the norm $\nrd(q)$
is positive at both of the roots, and antiholomorphically when it is negative at both places. Note that at a ramified real place (so for example at $t_0$ in our cubic example) the reduced norm $\nrd(q)>0$ for all nonzero $q$.

If $q\ne0$ and $\trd(q)=0$, then $q^2=-\nrd(q)$ is a nonzero scalar,
so its projective action is an involution.
If moreover $\nrd(q)$ is positive at both split places, $q$ acts as a half-turn (a rotation through $\pi$)
in each factor of $\HH\times\HH$ and fixes a point; if it is negative at both,
$q$ acts as a reflection in a geodesic in each factor and fixes a flat plane.
Section~5 uses such elements to compute fundamental groups.

An \emph{order} $O\subset\cB$ is a subring containing $1$ which is a finitely generated $R$-module spanning $\cB$
over $F$. It has rank four over $R$.
For example  $R+iR+jR+ijR$ is always an order when $a,b\in R$.
If $O$ is written in a basis, $e_1,...,e_4$,
$$O = e_1R+e_2R+e_3R+e_4R,$$
then one  needs to show that $1\in O$ and  $e_ie_j\in O$.
An order is \emph{maximal} if it is not
properly contained
in another order. Every order is contained in a maximal order.

For an order $O$ with an $R$-basis $e_1,\ldots,e_4$, define
$$ \delta(O)=\sqrt{\left|N_{F/\QQ}\left(\det(\trd(e_i e_j))\right)\right|}.$$
Here $N_{F/\QQ}$ is the field norm, computed by multiplying the values at all real embeddings.
 If $M$ is a maximal order containing
 $O$, then the additive index is
$[M:O]=\delta(O)/\delta(M)$. Note that $\delta(M)$ depends only on $\cB$, since it is equal to the
product of $|R/\mathfrak p|$ over the finitely
many primes $\mathfrak p$ at which
$\cB$ ramifies \cite[Chapter~15]{Voight}.

The norm-one elements form a group $O^1$, whose
projective image is
$$ \Lambda=O^1/\{\pm1\} \subset\PSL_2(\RR)\times\PSL_2(\RR).$$

By \cite{BorelHC,Shavel}, $\Lambda$ is an irreducible lattice in
$\PSL_2(\RR)\times\PSL_2(\RR)$, cocompact because $\cB$ is a division
algebra. The same holds for every group $\Gamma$ commensurable with $\Lambda$,
that is, with $\Gamma\cap\Lambda$ of finite index in both. So by
Proposition~\ref{lattice}, such a $\Gamma$ gives a fake quadric as soon as it
is torsion-free and $\vol(\Gamma\backslash(\HH\times\HH))=16\pi^2$.

\section{The quaternion algebra and the order}

In \cite{LSV}, Linowitz, Stover, and Voight  give a finite list of   arithmetic lattices that
can contain the fundamental group of an irreducible fake quadric. Some of these
are worked out explicitly, see \cite{Shavel,Takeuchi90,Dzambic}. Other entries may contain no
fake quadric, or more than one. For a general
entry the remaining work involves the study of torsion-free groups in the commensurability class of
the norm-one group of the order.
The algebra of this section belongs to the entry
$(n,d,D,N)=(3,169,13,1)$ in \cite{LSV}.

The field is
$$ F=\QQ(t),\qquad t^3-t^2-4t-1=0,$$
and there are
several useful maps. By substituting  a root $t_\nu$ of $f$ for $t$  we get  $\sigma_\nu\colon F\to\RR$.
Let $R=\ZZ[t]$. Since $f(-4)$ is divisible by $13$, substituting $t=-4$ and reducing
modulo $13$ gives a ring homomorphism $R\to\FF_{13}$. Let
$\mathfrak p=(13,t+4)$ denote the kernel of this map.

The quaternion algebra is
$$ \cB=(a,b)_F, \ \ a=3t^2-4t-12,\ \  b=-2t^2+4t+4.$$

At $t_0$ both $a$ and $b$ are negative, so the norm 
$x_0^2-ax_1^2-bx_2^2+abx_3^2$ is positive definite there. So in fact
$\cB$ is a division algebra. At $t_1$ and $t_2$ we have $a< 0 <b $. 

Define
\[
e_0=1,\qquad e_1=i,\qquad
e_2=\frac{t+i+j}{2},\qquad
e_3=\frac{-2ti+tj+ij}{4}.
\]
Then
\[
M=Re_0+Re_1+Re_2+Re_3
\]
is an order: it contains $1$, and direct multiplication shows that
$e_pe_q\in M$ for all $p,q$. A direct computation gives that $\delta(M) =13$. Since $\cB$ ramifies at exactly one
real place it also ramifies at an odd number of nonzero prime ideals of $R$. So we see that $\delta$ of a maximal order can't
be equal to $1$. It follows that $M$ is maximal.

The norm-one elements
in $M$, denoted as $M^1$ give a group.
Let
$$\Lambda = M^1/\{\pm 1\} \subset \PSL_2(\RR)\times \PSL_2(\RR)$$
be its projective image under  the embeddings $\sigma_1, \sigma_2$. By \cite{BorelHC,Shavel,LSV} it is an
irreducible lattice, cocompact since $\cB$ is a division algebra.

\begin{lemma}\label{notorsion}
  The group $\Lambda$ is torsion-free.
  \end{lemma}
  
\begin{proof} We want to rule out elements $q\neq  \pm 1$ in $M^1$ for $q^p=\pm 1$ for $p$ a prime.
  For odd primes we can assume $q^p=1$.

  Let $q=x_0+x_1i+x_2j+x_3ij$.
  If $q= x_0$ then $\nrd(q)=x_0^2=1$, contradicting $q \neq \pm 1$.
  We also have the identity
  $$q^2-\trd(q)q+\nrd(q)=0.$$
  
 So $\QQ(t)(q)$ is a field $L$ with $[L:\QQ]=6$.
For odd $p$,  the field $\QQ(q)\subset L$ has degree $p-1$. Since $p-1$ divides $6$, we get that  $p$ is either $3$ or $7$. In the case $p=7$
the two fields agree and $\QQ(t)$ would be isomorphic to $\QQ(z)$ with $z =e^{2\pi i/7}+ e^{-2\pi i/7}$,
however the two fields have different discriminants ($169$ and $49$ respectively). 

It remains to consider the cases $q^2=-1$ and $q^3=1$. 

Reduction modulo $\mathfrak p$ sends $t$ to $-4$, $a$ to $0$, and $b$ to $8$.
Since $8$ is a nonsquare in $\FF_{13}$, put
\[
\FF_{13^2}=\FF_{13}(w),\qquad w^2=8.
\]
The assignments $i\mapsto0$, $j\mapsto w$ respect the quaternion relations.
Since the denominators in the basis of $M$ are invertible modulo $13$,
they define a ring homomorphism
\[
H:M\longrightarrow\FF_{13^2}.
\]

Conjugation descends to $w\mapsto-w$. Thus, for $q\in M^1$ with
$H(q)=P+Qw$,
\[
1=H(q\overline q)=(P+Qw)(P-Qw)=P^2-8Q^2.
\]
This gives $14$ solutions, a group $\ZZ/14$ generated by $3+w$. If $q^2=-1$, then $H(q)^2=-1$, so $H(q)$ has order $4$.
If $q^3=1$ and $q\ne1$, the division-algebra property gives $q^2+q+1=0$.
Hence
\[
H(q)^2+H(q)+1=0,
\]
so $H(q)\ne1$ in characteristic $13$ and has order $3$.
Neither order is possible in a group of order $14$. This proves the lemma.
\end{proof}

The matrices in the introduction are the images of the following elements
of $M$, written in the $e_0,e_1,e_2,e_3$ basis:
\[
\begin{aligned}
A&=[t^2+t,t^2-3,4-2t^2,0],\\
B&=[t^2+2t,t+1,-3t-4,0],\\
C&=[-t-2,-t,t,-1].
\end{aligned}
\]
A direct computation shows that $A,B,C$ have reduced norm one, so they lie
in $M^1$ and $\Gamma\subset\Lambda$.
We will prove the equality in the next section.

\section{Volume Estimates}

We use the product hyperbolic volume form, with curvature $-1$ on each
factor. For distance comparisons we use the maximum metric
\[
d(z,x)=\max\{d_{\HH}(z_1,x_1),d_{\HH}(z_2,x_2)\},
\]
where $d_{\HH}$ is hyperbolic distance.

The following is straightforward.

\begin{proposition}\label{dirichlet}
Let $\Gamma\subset\PSL_2(\RR)\times\PSL_2(\RR)$ be discrete and torsion-free,
$p\in\HH\times \HH$, and $E\subset\Gamma$ finite. The set
$$D(p,E)=\{x\in\HH^2:\ d(p,x)\le d(gp,x)\ \text{for all}\ g\in E\}$$
is star-shaped with respect to $p$, hence path-connected and contractible. If it is bounded,
then $\Gamma\backslash\HH^2$ is compact and
$\vol(\Gamma\backslash\HH^2)\le\vol(D(p,E))$.
\end{proposition}

In practice, now working with $\Gamma =\langle A, B, C\rangle$ and $p=(i,i)$, we are looking for a finite set of words $E$ and a radius $R$ so that
the sphere $S(p,R)$ is disjoint from $D(p,E)$. This can be proved by subdividing $S(p,R)$ into finitely many small regions and showing that
each can be mapped inside a smaller ball around $p$ by some $e\in E$. This then proves compactness and a finite volume.

Let $S=\Lambda\backslash\HH^2$. Gauss-Bonnet, $K_S^2=2e(S)$, and
Noether's formula give
\[
\vol(S)=4\pi^2e(S)=16\pi^2\chi(\mathcal O_S).
\]
Thus $\vol(S)=16\pi^2m$ for a positive integer $m$.
Once $D(p,E)$ is bounded, $\Gamma$ has finite index $n$ in $\Lambda$, and
\[
\vol(\Gamma\backslash\HH^2)=16\pi^2mn.
\]
Consequently, an upper bound strictly below $32\pi^2$ forces $m=n=1$.
This proves both $\Gamma=\Lambda$ and the required covolume $16\pi^2$.

Subdivide $B(p,R)$ into regions $V_i$. If there is a $g\in E$ so that
$${\rm sup}_{x \in V_i} d(gp,x) < {\rm inf}_{x\in V_i} d(p,x),$$
then $V_i$ is disjoint from  $D(p,E)$.

Adding up
the volumes of the remaining regions gives an
upper bound for the volume of $D(p,E)$.  If the bound is too large we can add more elements to $E$ or subdivide
the radius $R$ ball further.

The accompanying program \texttt{volume.py} uses polar coordinate boxes  and verifies the bounds below using the $138$ words
in \texttt{words.txt}, all of length at most $17$. This then gives the following.

\begin{lemma}\label{lem:computation}
Let $p=(\ii,\ii)$, and let $E\subset\Gamma$ consist of all elements
represented by words of length at most $17$ in $A,B,C$ and their inverses.
Then
\[
D(p,E)\subset B(p,7/2),\qquad \vol(D(p,E))<245<32\pi^2.
\]
\end{lemma}

As discussed above the volume estimate finishes the proof of Theorem~\ref{thm1}.

Now let $X = \Gamma \backslash (\HH\times \HH)$. The group $\Gamma$ is generated by three elements $A,B,C$. Since these are given by matrices
and $\Gamma$ acts freely on $\HH\times \HH$ we can search for short relations. For example, the following relations hold in the projective group $\Gamma$:
\begin{align}
C^{-2}B^{-1}C^{-2}B^{-1}C^2AC^{-2}B^{-1}AB^{-2}&=1,
\label{eq:relation-one}\\
(C^{-2}BC^{-1}B^{-1}A)^2CA^2CA^{-1}&=1,
\label{eq:relation-two}\\
C^{-3}B^{-1}AC^{-1}B^{-1}CBCB^{-1}C^{-2}B^{-1}CACAC^{-1}A^{-1}&=1.
\label{eq:relation-three}
\end{align}

So in $H_1(X)$ we get
$$ [B]=4[A],\ \ [C]=-[A],\ \ 7[A]=0. $$
On the other hand the reduction map
of Section~3 induces a homomorphism
\[
\Gamma=M^1/\{\pm1\}\longrightarrow\ZZ/7.
\]
Since $H(A)=3-w\ne\pm1$, this homomorphism is nontrivial and hence
surjective. Therefore
\[
H_1(X;\ZZ)\cong\ZZ/7.
\] The remaining assertions of Theorem~\ref{thm2} follow from
Theorem~\ref{thm1} and Proposition~\ref{lattice}.

\section{Involutions and fundamental group calculations}

Conjugation by $i$ fixes $e_0=1$ and $e_1=i$, and satisfies
\[
ie_2i^{-1}=t+i-e_2,\qquad ie_3i^{-1}=-ti-e_3.
\]
It therefore preserves the maximal order $M$ and normalizes
$\Gamma=M^1/\{\pm1\}$. Since $i^2=a$ and $\nrd(i)=-a$ is positive
at both split real places, $i$ induces a holomorphic involution on $X$. Set
\[
\Gamma_2=\langle\Gamma,i\rangle.
\]
Then $X/i=\Gamma_2\backslash(\HH\times\HH)$.

In a matrix form at either split
embedding:
\[
\begin{pmatrix}
0&-\alpha\\
\alpha&0
\end{pmatrix}
=\alpha J,
\qquad
J=\begin{pmatrix}0&-1\\1&0\end{pmatrix}.
\]
Since scalar multiples induce the same projective action,
the involution is represented by the pair $(J,J)$.
Its action on $\HH\times\HH$ is
\[
(z_1,z_2)\longmapsto
\left(-\frac{1}{z_1},-\frac{1}{z_2}\right).
\]

This action has order two and  fixes the point $p=(\ii,\ii)$.

Since $ip=p$, we have $(Ci)p=Cp$.
Moreover, $(Ci)^2=1$ projectively, so $Ci$ interchanges $p$ and $Cp$.
It therefore fixes their midpoint, obtained by taking the hyperbolic
midpoint in each factor.

The existence of a fixed point follows also from a more general argument: $Ci$
has reduced trace zero and reduced norm $-a$ (that is positive at both split real places).
Thus it acts as a  rotation through $\pi$ in both factors
and has a fixed point.

By Armstrong's theorem \cite{Armstrong}, $\pi_1(X/i)$ is the quotient of
$\Gamma_2$ by the normal subgroup generated by elements having fixed points.
The two elements above give $i=1$ and $C=1$ in this quotient. With $C=1$,
relations \eqref{eq:relation-one}--\eqref{eq:relation-three} reduce to
\[
B^{-2}AB^{-1}AB^{-2}=1,\qquad A^3=1,\qquad B^{-1}AB^{-2}A=1.
\]

Using $P=B^{-1}A$ we get $B=P^2$ and $A=P^3$, 
$P^4=P^9=1$, showing that $A=B=C=1$ and $\pi_1(X/i)=1$.

Every nontrivial stabilizer acts as a rotation through $\pi$ in both factors, so the fixed
points are isolated and give ordinary double points in the quotient.
The involution $i$ preserves each of the two hyperbolic factors and its
area form. Since these cohomology classes form a basis of $H^2(X;\RR)$ it follows that $i$ acts trivially
on $H^2(X;\RR)$.
The Lefschetz fixed point theorem implies that  $X/i$ has four ordinary double points.

By \cite[Proposition~5.22]{DR}, the minimal resolution $Z$ is minimal of
general type with $p_g(Z)=0$ and $K_Z^2=4$. Resolving the four ordinary
double points introduces four disjoint $(-2)$-spheres and does not change
the fundamental group. This proves Theorem~\ref{thm3}.

For the antiholomorphic involution one can use
\[
T_0=(2+4t-2t^2)i+(3-t^2)j-t^2ij.
\]
Note that conjugation by $T_0$ preserves $M$.
Also  $T_0$ has reduced trace zero and
\[
\nrd(T_0)=16(t^2-3t-2),
\]
which is negative at both split real places. Since $T_0^2$ is scalar,
it therefore induces an antiholomorphic involution $\tau$ on $X$.
Set $Y=X/\langle\tau\rangle$.

The two elements
\[
T_0,\qquad CT_0
\]
have trace zero and the same norm as $T_0$. Each is a reflection in both
factors and has a fixed plane. Now using Armstrong's formula again we see that in $\pi_1(Y)$
$T_0=1$ and $C=1$. Using the relations for $\pi_1(X)$ we get 
$A=B=1$. It follows that
$\pi_1(Y)=1$.

We still need to look at the fixed point set. If a nonidentity
element in $\langle\Gamma,T_0\rangle$ fixes a point, then
it is a geodesic reflection in each factor.
Distinct reflections have disjoint fixed planes: otherwise
their product would be a nonidentity element of $\Gamma$
fixing a point.

Each component of $\mathcal F=\Fix(\tau)$ is a compact quotient
of a plane $L_1\times L_2$ by a subgroup of $\Gamma$ preserving
both geodesics. Let $q\in M^1$ represent an
element in this  subgroup.
If its action reverses one of the geodesics, then its
matrix at that embedding has trace zero. 
However then its reduced
trace $\trd(q)$ also vanishes in $F$. It follows that $q^2=- 1$. Its projective class
would have order two and that contradicts $\Gamma$ being torsion free.
So the subgroup acts by translations and the quotient is a flat torus.

Since $\mathcal F$ is a union of tori, we get $e(Y)=2$.
Since  $\pi_1(Y)=1$, $Y$ is a homotopy $4$-sphere and
Freedman's theorem \cite{Freedman} implies that it is homeomorphic
to $S^4$.

Finally we can use Smith's fixed point theorem. We have a $\ZZ/2$ action. Since $|H_1(X;\ZZ)|$ is odd, and $e(X)=4$ we see that
the sum of the mod $2$ Betti numbers of $X$ is $4$. This shows that we can't have more than one torus as the branch locus.
This proves Theorem~\ref{thm4}.


\begin{thebibliography}{99}
\bibitem{Armstrong} M. A. Armstrong, \emph{The fundamental group of the orbit space of a discontinuous group}, Proc. Cambridge Philos. Soc. \textbf{64} (1968), 299-301.


\bibitem{Arnold} V. I. Arnold, \emph{A branched covering $\mathbb{CP}^2\to S^4$, hyperbolicity and projective topology}, Siberian Math. J. \textbf{29} (1988), 717-726.


 \bibitem{Cat2} I. Bauer, F. Catanese, F. Grunewald, \emph{The classification of surfaces with $p_g=q=0$ isogenous to a product of curves}, Pure Appl. Math. Q. \textbf{4} (2008), 547-586.


 \bibitem{BCP} I. Bauer, F. Catanese, R. Pignatelli, \emph{Surfaces of general type with geometric genus zero: a survey}, in: Complex and Differential Geometry, Springer Proc. Math. \textbf{8}, Springer, Berlin, 2011, 1-48.



\bibitem{BorelHC} A. Borel, Harish-Chandra, \emph{Arithmetic subgroups of algebraic groups}, Ann. of Math. (2) \textbf{75} (1962), 485-535.

\bibitem{BR} L. Borisov, C. Rito, \emph{Explicit Bicanonical Models of Eight Fake Quadrics}, arXiv:2608.12296 (2026).



\bibitem{Dzambic} A. Dzambic, \emph{Fake quadrics from irreducible lattices acting on the product of upper half planes}, Math. Ann. \textbf{360} (2014), 23-51.


\bibitem{DR} A. Dzambic, X. Roulleau, \emph{Automorphisms and quotients of quaternionic fake quadrics}, Pacific J. Math. \textbf{267} (2014), 91-120.

\bibitem{Finashin} S. Finashin, \emph{Exotic embeddings of $\# 6\mathbb{RP}^2$ in the $4$-sphere}, in: Proceedings of G\"okova Geometry-Topology Conference 2008, G\"okova, 2009, 151-169.

\bibitem{FKV2} S. M. Finashin, M. Kreck, O. Ya. Viro, \emph{Nondiffeomorphic but homeomorphic knottings of surfaces in the $4$-sphere}, in: Topology and geometry---Rohlin Seminar, Lecture Notes in Math. \textbf{1346}, Springer, 1988, 157-198.

\bibitem{Freedman} M. H. Freedman, \emph{The topology of four-dimensional manifolds}, J. Differential Geom. \textbf{17} (1982), 357-453.

\bibitem{FM} R. Friedman, J. W. Morgan, \emph{Algebraic surfaces and Seiberg-Witten invariants}, J. Algebraic Geom. \textbf{6} (1997), 445-479.



\bibitem{LP} Y. Lee, J. Park, \emph{A simply connected surface of general type with $p_g=0$ and $K^2=2$}, Invent. Math. \textbf{170} (2007), 483-505.

\bibitem{LSV} B. Linowitz, M. Stover, J. Voight, \emph{Commensurability classes of fake quadrics}, Selecta Math. (N.S.) \textbf{25} (2019), no. 3, Paper No. 48; Correction, ibid. \textbf{25} (2019), no. 4.



\bibitem{MORSU} G. Mati\'c, F. \"Ozt\"urk, J. Reyes, A. I. Stipsicz, G. Urz\'ua, \emph{An exotic $5\mathbb{RP}^2$ in the $4$-sphere}, arXiv:2312.03617.

\bibitem{PPS4} H. Park, J. Park, D. Shin, \emph{A simply connected surface of general type with $p_g=0$ and $K^2=4$}, Geom. Topol. \textbf{13} (2009), 1483-1494.

\bibitem{Rito} C. Rito, \emph{Computation of Singular Godeaux Surfaces and a New Explicit Fake Quadric}, with an appendix by C. Gleissner and N. Ruhland, arXiv:2509.08198 (2025).

\bibitem{Shavel} I. H. Shavel, \emph{A class of algebraic surfaces of general type constructed from quaternion algebras}, Pacific J. Math. \textbf{76} (1978), 221-245.


\bibitem{SS} A. I. Stipsicz, Z. Szab\'o, \emph{The real quotient of the Cartwright-Steger surface}, arXiv:2609.09207.

\bibitem{Takeuchi90} K. Takeuchi, \emph{Algebraic surfaces derived from unit groups of quaternion algebras}, in: Number Theory (Banff, AB, 1988), de Gruyter, Berlin, 1990, 529-552.

\bibitem{Voight} J. Voight, \emph{Quaternion algebras}, Graduate Texts in Math. \textbf{288}, Springer, 2021.

\end{thebibliography}
\end{document}